\documentclass[prd,showpacs,twocolumn]{revtex4-1}

\usepackage{amsmath}
\usepackage{amsfonts}
\usepackage{amsthm}
\usepackage{graphicx}
\usepackage{url}

\newtheorem{lemma}{Lemma}
\newtheorem{corollary}{Corollary}

\newcommand{\RR}{{\mathbb R}}
\newcommand{\CC}{{\mathbb C}}
\newcommand{\HH}{{\mathbb H}}

\newcommand{\OO}{{\mathbb O}}

\newcommand{\isom}{\cong}
\newcommand{\tr}{\mathrm{tr}}

\newcommand{\so}{\mathfrak{so}}
\newcommand{\su}{\mathfrak{su}}
\renewcommand{\sl}{\mathfrak{sl}}
\renewcommand{\sp}{\mathfrak{sp}}
\renewcommand{\aa}{\mathfrak{a}}
\newcommand{\cc}{\mathfrak{c}}
\newcommand{\dd}{\mathfrak{d}}
\newcommand{\ee}{\mathfrak{e}}
\newcommand{\ff}{\mathfrak{f}}
\renewcommand{\gg}{\mathfrak{g}}	
\newcommand{\der}{\mathfrak{der}}

\newcommand{\EE}{E}

\newcommand{\NN}{\mathcal{N}}
\newcommand{\UU}{\mathcal{U}}

\newcommand{\eA}{A}	
\newcommand{\eG}{G}	
\newcommand{\eS}{S}	
\newcommand{\eD}{D}	
\newcommand{\eE}{E}	
\newcommand{\eF}{F}	
\newcommand{\eX}{X}	
\newcommand{\eY}{Y}	
\newcommand{\eZ}{Z}	

\newcommand{\HHH}{\mathrm{H}_3(\OO)}
\newcommand{\AAA}{{\cal A}}
\newcommand{\BBB}{{\cal B}}
\newcommand{\EEE}{{\cal E}}
\newcommand{\III}{{\cal I}}
\newcommand{\PPP}{{\cal P}}
\newcommand{\XXX}{{\cal X}}
\newcommand{\YYY}{{\cal Y}}
\newcommand{\Pc}{\hbox{\boldmath${\cal P}$}}

\renewcommand{\Im}{\mathrm{Im}\,}
\renewcommand{\bar}[1]{\overline{#1}}
\newcommand{\mb}[1]{\mathbf{#1}}

\renewcommand{\footnote}[1]{(#1)}

\usepackage{xcolor}
  \definecolor{bluegray}{rgb}{0.4,0.6,0.8}
  \definecolor{forest}{rgb}{0,0.5,0}

\begin{document}


\title{\boldmath A New Division Algebra Representation of $\EE_7$}

\author{Tevian Dray}
\email{tevian@math.oregonstate.edu}
\affiliation{Department of Mathematics, Oregon State University,
Corvallis, OR  97331,USA}

\author{Corinne A. Manogue}
\email{corinne@physics.oregonstate.edu}
\affiliation{Department of Physics, Oregon State University,
Corvallis, OR  97331, USA}

\author{Robert A. Wilson}
\email{r.a.wilson@qmul.ac.uk}
\affiliation{School of Mathematical Sciences, Queen Mary,
University of London, London E1 4NS, UK}

\begin{abstract}
We decompose the Lie algebra $\ee_{8(-24)}$ into representations of
$\ee_{7(-25)}\oplus\sl(2,\RR)$ using our recent description of $\ee_8$ in
terms of (generalized) $3\times3$ matrices over pairs of division algebras.
Freudenthal's description of both $\ee_7$ and its minimal representation are
therefore realized explicitly within $\ee_8$, with the action given by the
(generalized) matrix commutator in $\ee_8$, and with a natural
parameterization using division algebras.  Along the way, we show how to
implement standard operations on the Albert algebra such as trace of the
Jordan product, the Freudenthal product, and the determinant, all using
commutators in $\ee_8$.
\end{abstract}

\maketitle

\section{Introduction}
\label{intro}
The Freudenthal--Tits magic square~\cite{Freudenthal,Tits} parameterizes
certain Lie algebras, including most of the exceptional Lie algebras, in terms
of two division algebras.  Building on previous work~\cite{Sudbery,%
SudberyBarton,Lorentz,Denver,York,AaronThesis,Sub,Structure,Cubies,%
JoshuaThesis,SO42,2x2}, we recently provided~\cite{Magic,Octions} a
description of $\ee_8$---and hence the entire magic square---in terms of
$3\times3$ (generalized) matrices over the tensor product of two copies of the
octonions $\OO$ (or their split cousins $\OO'$).  This construction, which is
summarized in Section~\ref{e8}, provides an explicit implementation of the
Vinberg description~\cite{Vinberg} of the magic square, the ``half-split''
form of which is shown in Table~\ref{3x3}.

\begin{table*}
\small
\begin{center}
\begin{tabular}{|c|c|c|c|c|}
\hline
&$\RR$&$\CC$&$\HH$&$\OO$\\\hline
$\RR'$
&$\su(3,\RR)$&$\su(3,\CC)$&$\cc_3\isom\su(3,\HH)$&$\ff_4\isom\su(3,\OO)$\\
\hline
$\CC'$&
$\sl(3,\RR)$&$\sl(3,\CC)$&$\aa_{5(-7)}\isom\sl(3,\HH)$
  &$\ee_{6(-26)}\isom\sl(3,\OO)$\\
\hline
$\HH'$&
$\cc_{3(3)}\isom\sp(6,\RR)$&$\su(3,3,\CC)\isom\sp(6,\CC)$
  &$\dd_{6(-6)}\isom\sp(6,\HH)$&$\ee_{7(-25)}\isom\sp(6,\OO)$\\
\hline
$\OO'$&
$\ff_{4(4)}$&$\ee_{6(2)}$&$\ee_{7(-5)}$&$\ee_{8(-24)}$\\
\hline
\end{tabular}
\end{center}
\caption{The ``half-split'' real form of the Freudenthal--Tits magic square of
Lie algebras.}
\label{3x3}
\end{table*}

Since the minimal representation of $\ee_8$ is the adjoint representation, the
treatment in~\cite{Magic,Octions} therefore provide a matrix description of
the adjoint representation of~$\ee_8$, which extends naturally to the Lie
group $\EE_8$.
Working backward, it is straightforward to reduce this description to the
adjoint representation of each of the other Lie algebras in the magic square.
But we can do more: Breaking the adjoint representation of one algebra down to
the adjoint representation of a subalgebra results in a decomposition of the
original algebra into several pieces, and it is instructive to examine these
pieces explicitly.

In this paper, we apply this idea to the well-known decomposition of $\ee_8$
over $\ee_7\oplus\su(2)$, including a discussion of different possible real
forms, corresponding to the use of either split or regular division algebras.
(In a separate paper~\cite{MagicE6}, we provided a similar discussion of
the decomposition of $\ee_8$ over $\ee_6\oplus\su(3)$.)
We thus obtain a description of both $\ee_7$ and its minimal representation
in terms of $\ee_8$, which we then compare with Freudenthal's
description~\cite{FreudenthalE7}.  Since this latter description involves
operations on elements of the Albert algebra, including the Freudenthal
product, the trace of the Jordan product, and the determinant, we obtain
explicit implementations of these operations as commutators in $\ee_8$.

In addition to its mathematical interest, this construction may be useful to
those who are trying to build the Standard Model of particle physics using the
exceptional Lie groups; see for
instance~\cite{Gunaydin,GurseyE6,Bars80,Lisi,Chester,Magic} and references
cited there.
%
%

After briefly reviewing our description of $\ee_8$ in Section~\ref{e8}, we
summarize Freudenthal's description of $\ee_7$ in Section~\ref{e7Freud}.  In
Section~\ref{e7}, we put these pieces together, showing how the objects and
operations in Freudenthal's description of the minimal representation of
$\ee_7$ can be implemented in terms of commutators in $\ee_8$.  In
Section~\ref{e7adj}, we briefly return to the adjoint representation of
$\ee_7$ as a subalgebra of $\ee_8$, showing how to implement Freudenthal's
description in this setting.  Finally, in Section~\ref{discussion}, we discuss
our results and suggest some possible applications to the physics of the
Standard Model.

\section{\boldmath Adjoint representation of $\ee_8$}
\label{e8}


\begin{figure}
\centering
\includegraphics[width=5cm]{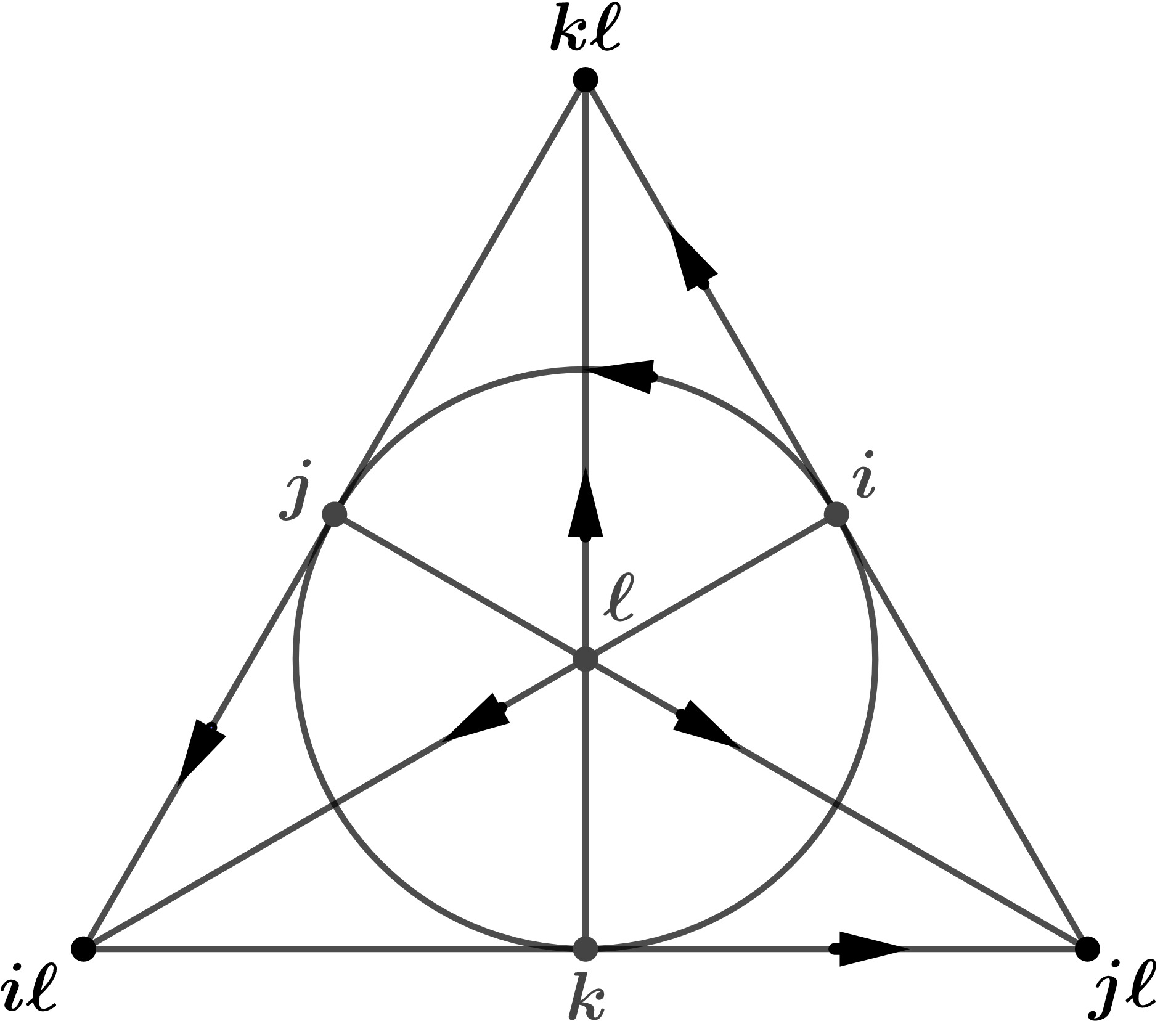}
\caption{A graphical representation of the  octonionic multiplication table.}
\label{omult3}
\end{figure}

\begin{table}
\centering
\small
\begin{tabular}[b]{|c|c|c|c|c|c|c|c|}
\hline
&\boldmath$i$&\boldmath$j$&\boldmath$k$&\boldmath$k\ell$
  &\boldmath$j\ell$&\boldmath$i\ell$&\boldmath$\ell$\\\hline
\boldmath$i$&$-1$&$k$&$-j$&$j\ell$&$-k\ell$&$\ell$&$i\ell$\\\hline
\boldmath$j$&$-k$&$-1$&$i$&$-i\ell$&$\ell$&$k\ell$&$j\ell$\\\hline
\boldmath$k$&$j$&$-i$&$-1$&$-\ell$&$i\ell$&$-j\ell$&$k\ell$\\\hline
\boldmath$k\ell$&$-j\ell$&$i\ell$&$\ell$&$-1$&$i$&$-j$&$-k$\\\hline
\boldmath$j\ell$&$k\ell$&$\ell$&$-i\ell$&$-i$&$-1$&$k$&$-j$\\\hline
\boldmath$i\ell$&$\ell$&$-k\ell$&$j\ell$&$j$&$-k$&$-1$&$-i$\\\hline
\boldmath$\ell$&$-i\ell$&$-j\ell$&$-k\ell$&$k$&$j$&$i$&$-1$\\\hline
\noalign{\vspace{0.15in}}
\end{tabular}
\caption{The octonionic multiplication table.}
\label{omult}
\end{table}

The \textit{octonions} $\OO$ are the real algebra spanned by the identity
element~$1$ and seven square roots of~$-1$ that we denote
$\UU=\{i,j,k$,$k\ell,j\ell,i\ell,\ell\}$, whose multiplication table is neatly
described by the oriented Fano geometry shown in Figure~\ref{omult3}, and
given explicitly in Table~\ref{omult}.  The \textit{split octonions} $\OO'$
are the real algebra spanned by the identity element $1$, together with seven
square roots of~$\pm1$ that we denote $\UU'=\{I,J,K,KL,JL,IL,L\}$, whose
multiplication table is given in Table~\ref{smult}.  Our use of $1$ for the
identity elements of both $\OO$ and $\OO'$ should be clear from the context.

\begin{table}
\centering
\small
\begin{tabular}{|c|c|c|c|c|c|c|c|}
\hline
$$&\boldmath$I$&\boldmath$J$&\boldmath$K$
  &\boldmath$KL$&\boldmath$JL$&\boldmath$IL$&\boldmath$L$\\\hline
\boldmath$I$&$-1$&$K$&$-J$&$JL$&$-KL$&$-L$&$IL$\\\hline
\boldmath$J$&$-K$&$-1$&$I$&$-IL$&$-L$&$KL$&$JL$\\\hline
\boldmath$K$&$J$&$-I$&$-1$&$-L$&$IL$&$-JL$&$KL$\\\hline
\boldmath$KL$&$-JL$&$IL$&$L$&$1$&$-I$&$J$&$K$\\\hline
\boldmath$JL$&$KL$&$L$&$-IL$&$I$&$1$&$-K$&$J$\\\hline
\boldmath$IL$&$L$&$-KL$&$JL$&$-J$&$K$&$1$&$I$\\\hline
\boldmath$L$&$-IL$&$-JL$&$-KL$&$-K$&$-J$&$-I$&$1$\\\hline
\end{tabular}
\caption{The split octonionic multiplication table.}
\label{smult}
\end{table}



Our matrix representation of $\ee_8$ was given in~\cite{Magic} and further
discussed in~\cite{MagicE6}; here we present only a brief summary in order to
introduce our notation.

The 248 elements of the adjoint representation of $\ee_{8(-24)}$ can be
represented as the \hbox{$3\times8\times8+2\times14=192+28=220$} tracefree
$3\times3$ anti-Hermitian matrices over $\OO'\otimes\OO$, together with the
$14+14=28$ elements of $\gg_{2(2)}\oplus\gg_2$.
These tracefree matrices take the form
\begin{align}
\eX_p &= \begin{pmatrix}0& p& 0\\ -\bar{p}& 0 &0\\ 0& 0& 0\\\end{pmatrix} ,
\qquad
\eD_q = \begin{pmatrix}q& 0& 0\\ 0& -q& 0\\ 0& 0& 0\\\end{pmatrix} ,\\
\displaybreak[0]
\eY_p &= \begin{pmatrix}0& 0& 0\\ 0& 0& p\\ 0& -\bar{p} &0\\\end{pmatrix} ,
\qquad
\eE_q = \begin{pmatrix}0& 0& 0\\ 0& q& 0\\ 0& 0& -q\\\end{pmatrix} ,\\
\displaybreak[0]
\eZ_p &= \begin{pmatrix}0& 0& -\bar{p}\\  0& 0& 0\\ p& 0 &0\\\end{pmatrix} ,
\qquad
\eF_q = \begin{pmatrix}-q& 0& 0\\ 0& 0& 0\\ 0& 0& q\\\end{pmatrix} ,
\end{align}
with $p\in\OO'\otimes\OO$ and $q\in\Im\OO\cup\Im\OO'$.  We will usually
restrict our labels to elements of $\UU$ and $\UU'$, which suffice to generate
a basis of $\ee_8$.

As pointed out in~\cite{Lorentz}, the use of octonionic matrices to represent
the \textit{action} of a Lie (group or) algebra requires \textit{nesting}; the
action can not always be represented in terms of a single matrix.  The key to
the construction in~\cite{Magic} is that nested elements such as
\begin{equation}
\eD_{p,q} = \frac12 \> [\eD_p,\eD_q] = \frac12 \> [\eX_p,\eX_q]
\end{equation}
with $p\ne q\in\UU$ (or $p\ne q\in\UU'$) can nonetheless be represented as
\textit{generalized} $3\times3$ matrices, leading to a (generalized) matrix
representation of the adjoint representation.  This reinterpretation is
accomplished by writing
\begin{equation}
\eD_{p,q}
  = \begin{pmatrix}p\circ q& 0& 0\\ 0& p\circ q& 0\\ 0& 0& 0\\\end{pmatrix}
\label{comp}
\end{equation}
where $\circ$ denotes \textit{composition}, that is
\begin{equation}
(p\circ q)r = p(qr)
\end{equation}
for $p,q,r\in\OO$ (or $p,q,r\in\OO'$).  Composition is associative by
definition!  Furthermore, composition respects anticommutativity, in the sense
that
\begin{equation}
qp = -qp \Longrightarrow q\circ p=-p\circ q
\end{equation}
which holds in particular for $p\ne q\in\UU$ (or $p\ne q\in\UU'$).  For such
$p$, $q$,
\begin{align}
p\circ(p\circ q) &= p^2 q = -(p\circ q)\circ p \label{ppq} ,\\
(p\circ q)\circ(q\circ r) &= p\circ q^2 \circ r = q^2 p\circ r \label{pqr}
\end{align}
for any $r\in\OO$ (or $r\in\OO'$), since $q^2\in\RR$.
Similarly, defining
\begin{align}
\eE_{p,q} &= \frac12 \> [\eE_p,\eE_q] = \frac12 \> [\eY_p,\eY_q] ,\\
\eF_{p,q} &= \frac12 \> [\eF_p,\eF_q] = \frac12 \> [\eZ_p,\eZ_q] .
\end{align}
leads to
\begin{equation}
\eE_{p,q}
  = \begin{pmatrix}0& 0& 0\\ 0& p\circ q& 0\\ 0& 0& p\circ q\\\end{pmatrix} ,
\>
\eF_{p,q}
  = \begin{pmatrix}p\circ q& 0& 0\\ 0& 0& 0\\ 0& 0& p\circ q\\\end{pmatrix} ,
\end{equation}
again with $p\ne q\in\UU$ (or $p\ne q\in\UU'$).

We seem to have too many generators, but triality ensures that only $28+28=56$
of the $\eD$s, $\eE$s, and $\eF$s are independent.  We refer to the $\eD$s and
$\eX$s as being of \textit{type~I}, the $\eE$s and $\eY$s as \textit{type~II},
and the $\eF$s and $\eZ$s as \textit{type~III}.  We need all $64+64+64=192$ of
the $\eX$s, $\eY$s, and $\eZ$s, but it is enough to take the $28+28=56$
elements of a single type, such as the $\eD$s, to generate
$\so(4,4)\oplus\so(8)\subset\ee_8$.

As discussed in more detail in~\cite{MagicE6}, it is sometimes convenient to
replace the $21+21=42$ nested $\eD$s, generating (type~I)
$\so(3,4)\oplus\so(7)$, by a second, slightly different basis.  Of the 21
elements
\begin{align}
\eD_{p,q} + \eE_{p,q} + \eF_{p,q}
  &= 2\begin{pmatrix}
	p\circ q& 0& 0\\ 0& p\circ q& 0\\ 0& 0& p\circ q\\
  \end{pmatrix} \nonumber\\
  &= 2(p\circ q)\,\III
\label{DEF}
\end{align}
with $p\ne q\in\UU$, only 14 are independent.  In fact, for each element
$r\in\UU$, there are three pairs $p,q\in\UU$ such that $pq=r$, but the
corresponding elements of the form~(\ref{DEF}) span a two-dimensional vector
space, every element of which is a (nested) multiple of $\III$.  We choose a
particular orthogonal basis $\{\eA_r,\eG_r\}$ of each such space; these 14
elements generate $\gg_2$.  We choose the 7 elements $\eA_r$ so that they fix
$\ell$ (at the group level); the 8 $\eA$s, together with $\eG_\ell$, therefore
generate the $\su(3)\subset\gg_2$ that fixes $\ell\in\OO$.
\footnote{Both $\eA_\ell$ and $\eG_\ell$ fix $\ell$ at the group level, so
their exact form is conventional.}
A similar construction holds for $\OO'$.
We need 7 more elements to generate (type~I) $\so(7)$; triality allows us to
choose them in the unnested form
\begin{equation}
\eS_p
  = \eE_p - \eF_p
  = \begin{pmatrix}p& 0& 0\\ 0& p& 0\\ 0& 0& -2p\\\end{pmatrix} .
\end{equation}

Our second basis for $\ee_8$ therefore consists of 192 off-diagonal elements,
namely the $\eX$s, $\eY$s, and~$\eZ$s, together with the $4\times(7+7)=56$
(single-index) $\eD$s, $\eS$s, $\eG$s, and $\eA$s.  Only the $\eG$s and $\eA$s
are nested; both are (nested) multiples of the identity element.  In this
representation, one possible choice of generators of the Cartan subalgebra is
then $\{\eA_\ell,\eG_\ell,\eS_\ell,\eD_\ell\}$ and
$\{\eA_L,\eG_L,\eS_L,\eD_L\}$, corresponding to the $\so(8)$ and $\so(4,4)$
subalgebras of $\ee_8$, over $\OO$ and $\OO'$, respectively.

\section{\boldmath Freudenthal Description of $\ee_7$}
\label{e7Freud}

We summarize here our treatment~\cite{Cubies} of the description of
$\ee_7=\ee_{7(-25)}$ originally given by Freudenthal~\cite{FreudenthalE7}.

Let $\HHH$ denote the Albert algebra of $3\times3$ Hermitian matrices over
$\OO$.  The \textit{Jordan product} of $\XXX,\YYY\in\HHH$ is given by
\begin{equation}
\XXX\circ\YYY = \frac12 \Bigl( \XXX\YYY + \YYY\XXX \Bigr)
\end{equation}
and their \textit{Freudenthal product} is
\begin{align}
\XXX*\YYY
  = \XXX\circ\YYY
	&- \frac12\Bigl( (\tr\XXX)\,\YYY + (\tr\YYY)\,\XXX \Bigr)
  \nonumber\\
	&+ \frac12\Bigl( (\tr\XXX)(\tr\YYY)-\tr(\XXX\circ\YYY) \Bigr) \>\III
\end{align}
where $\III$ denotes the identity matrix.
\footnote{The use of $\circ$ to denote both composition and the Jordan product
should be clear from context.}

The generators of $\ee_6=\ee_{6(-26)}$ fall into one of three categories;
there are 26 \textit{boosts}, 14 \textit{derivations} (of $\OO$, that is,
elements of $\gg_2$), and 38 remaining \textit{rotations} (the remaining
generators of $\ff_4$).  For both boosts and rotations, $\phi\in\ee_6$ can be
treated as a $3\times3$, tracefree, octonionic matrix; boosts are Hermitian,
and rotations are anti-Hermitian.  Such matrices $\phi\in\ee_6$ act on the
Albert algebra via
\begin{equation}
\XXX \longmapsto \phi\XXX + \XXX\phi^\dagger
\label{e6act}
\end{equation}
where $\dagger$ denotes conjugate transpose (in $\OO$).  Derivations can be
obtained by successive rotations (or boosts) through \textit{nesting},
corresponding to commutators in the Lie algebra; it therefore suffices to
consider the boosts and rotations, that is, to consider matrix
transformations.
The dual representation of $\ee_6$ is formed by setting
\begin{equation}
\phi' = -\phi^\dagger
\label{phiadj}
\end{equation}
for both boosts and rotations, which leaves rotations alone but reverses the
sign of all boosts.

Freudenthal~\cite{FreudenthalE7} describes elements of $\ee_7$ as
\begin{equation}
\Theta = (\phi,\rho,\AAA,\BBB)
\label{ThetaDef}
\end{equation}
where $\phi\in\ee_6$, $\rho\in\RR$, and $\AAA,\BBB\in\HHH$.  Thus, $\ee_7$ is
spanned by the 78 elements of~$\ee_6$ ($\phi$), together with $27+27=54$
(null) translations ($\AAA,\BBB$), and a dilation ($\rho$); $\EE_7$ is
the conformal group of $\EE_6$.

What does $\Theta$ act on?  Freudenthal~\cite{FreudenthalE7} explicitly
constructs the minimal representation of~$\ee_7$, which consists of elements
of the form
\begin{equation}
\Pc=(\XXX,\YYY,p,q)
\label{Pcdef}
\end{equation}
where $\XXX,\YYY\in\HHH$, and $p,q\in\RR$.
%
Freudenthal also tells us that $\Theta$ acts on $\Pc$ via
\begin{align}
  \XXX &\longmapsto \phi(\XXX) + \frac13\,\rho\,\XXX + 2 \BBB*\YYY + \AAA\,q
, \label{FreudX}\\
  \YYY &\longmapsto 2 \AAA*\XXX + \phi'(\YYY) - \frac13\,\rho\,\YYY + \BBB\,p
, \label{FreudY}\\
  p &\longmapsto \tr(\AAA\circ \YYY) - \rho\,p
, \label{pFreud}\\
  q &\longmapsto \tr(\BBB\circ \XXX) + \rho\,q 
. \label{qFreud}
\end{align}
How are we to visualize this action?  Consider a fixed element
$\Theta=(0,0,\AAA,0)$ acting repeatedly on $\Pc=(0,0,0,q)$.  We have
\begin{align}
(0,0,\AAA,0): (0,0,0,q)
	&\longmapsto (q\,\AAA,0,0,0) \nonumber\\
	&\longmapsto (0,2q\,\AAA*\AAA,0,0) \nonumber\\
	&\longmapsto (0,0,6q\,\det\AAA,0) \nonumber\\
	&\longmapsto (0,0,0,0)
\label{Araise}
\end{align}
where we have used the fact that
\begin{equation}
\det(\AAA) = \frac13 \, \tr \Big( (\AAA*\AAA) \circ \AAA \Big) .
\label{Det}
\end{equation}
Similarly, we have
\begin{align}
(0,0,0,\BBB): (0,0,p,0)
	&\longmapsto (0,p\,\BBB,0,0) \nonumber\\
	&\longmapsto (2p\,\BBB*\BBB,0,0,0) \nonumber\\
	&\longmapsto (0,0,0,6p\,\det\BBB) \nonumber\\
	&\longmapsto (0,0,0,0) .
\label{Braise}
\end{align}
Thus, the null translations $\AAA$ and $\BBB$ act as raising and lowering
operators, and we therefore think of $\Pc$ as a ``tower''
which we henceforth refer to as a \textit{Freudenthal tower}, with
\textit{anchors} $p$ and $q$, as shown schematically in Table~\ref{tower}.

\begin{table}
\centering
\fbox{
\begin{tabular}{ccc}
& $p$ \\
& $|$ \\
$\uparrow$\quad\null & $\YYY$ \\
$\AAA$\quad\null & $|$ & \quad$\BBB$ \\
& $\XXX$ & \quad$\downarrow$ \\
& $|$ \\
& $q$
\end{tabular}
}
\caption{A schematic representation of elements $\Pc=(\XXX,\YYY,p,q)$ of the
minimal representation of $\ee_7=\ee_{7(-25)}$ as a \textit{Freudenthal
tower}, showing also the action of the null translations $(0,0,\AAA,0)$ and
$(0,0,0,\BBB)$.}
\label{tower}
\end{table}

\section{\boldmath A Representation of $\ee_7$ in $\ee_8$}
\label{e7}

\subsection{Overview}
\label{overview}

The adjoint representation of $\ee_8=\ee_{8(-24)}$ was given explicitly in
Section~\ref{e8} in terms of (generalized) $3\times3$ matrices over
$\OO'\otimes\OO$.  The adjoint representation of $\ee_7=\ee_{7(-25)}$ is
immediately obtained by restricting this construction to $\HH'\otimes\OO$.
Explicitly, we now have $3\times4\times8=96$ off-diagonal elements ($\eX_p$,
$\eY_p$, and $\eZ_p$, with $p\in\HH'\otimes\OO$), the $28$ diagonal elements
of $\so(8)$ ($\eA_q$, $\eG_q$, $\eS_q$, and $\eD_q$, with $q\in\Im\OO$), the 6
independent tracefree diagonal elements over $\HH'$ ($\eD_P$ and $\eS_P$, with
$P\in\Im\HH'$), and finally the 3 elements of $\der(\HH')$ ($(P\circ Q)\III$,
with $P,Q\in\Im\HH'$; these are the $\eG_P$, for $P\in\Im\HH'$), for a total
of $96+28+6+3=133$ elements in adjoint $\ee_7$, as expected.
In $\ee_8$, we also have the 3 elements of the $\sl(2,\RR)$ that centralizes
$\ee_7$; these are the $\eA_P$, with $P\in\Im\HH'$.

Alternatively, we can start by enumerating the basis elements of
$\ee_6=\ee_{6(-26)}$.  We take $\HH'$ to be spanned by $\{1,K,KL,L\}$, and
$\CC'$ to be the split complex subalgebra containing $L$.  We have
$3\times2\times8=48$ off-diagonal matrices over $\CC'\otimes\OO$, together
with the 28 diagonal elements of $\so(8)$ given above, together with the 2
independent tracefree diagonal elements over $\CC'$, namely $\eD_L$ and
$\eS_L$, for a total of $48+28+2=78$ elements.  The remaining 55 basis
elements of $\ee_7$ are the $2\times26=52$ tracefree matrices having $K$ or
$KL$ as a factor, along with $\eG_K$, $\eG_{KL}$, and $\eG_L$.

\subsection{\boldmath Freudenthal towers in $\ee_8$}

What elements are left in $\ee_8$?

We still have all elements labeled by $Q\in\{I,J,IL,JL\}$.  For each such
$Q$, there are 26 tracefree elements of $\ee_8$ with $Q$ as a factor, along
with $\eG_Q$ and $\eA_Q$, for a total of $4\times(26+1+1)=112$ elements of
$\ee_8$ that are neither in adjoint $\ee_7$ nor in adjoint $\sl(2,\RR)$---just
enough degrees of freedom to build \textit{two} Freudenthal towers.  We now
show how to do exactly that, thus decomposing $\ee_{8(-24)}$ explicitly over
$\ee_{7(-25)}\oplus\sl(2,\RR)$ as
\begin{equation}
\label{decomp}
\mb{248} = \mb{133} + 2\times\mb{56} + \mb{3} .
\end{equation}

As discussed in~\cite{MagicE6}, for each $Q\in\Im\OO'$, the 26 independent
tracefree elements having $Q$ as a factor, together with $\eG_Q$, span a copy
of the Albert algebra, with $-\frac12\eG_Q$ playing the role of ``$Q\,\III$''.
Thus, $\ee_{7(-25)}$ contains $\ee_{6(-26)}$ and two copies of the Albert
algebra (labeled by $K$ and $KL$), with a single element left over to be the
dilation, namely $\eG_L$.  Similarly, $\ee_{8(-24)}$ contains $\ee_{7(-25)}$,
$\sl(3,\RR)$, and 4 copies of the Albert algebra (labeled by
$Q\in\{I,J,IL,JL\}$), with 4 elements left over (the $\eA_Q$).

However, as with the representations of $\ee_6$ constructed in~\cite{MagicE6},
there's a subtlety in how we build these copies of the Albert algebra: We must
take $Q$ to be null.  Thus, the null translations of $\ee_{7(-25)}$ ($\AAA$ and
$\BBB$ of Section~\ref{e7Freud}) are labeled not by $K$ and $KL$, but rather
by $K\pm KL$.  Similarly, the Freudenthal towers on which $\ee_7$ acts must be
labeled by $\{I\pm IL,J\pm JL\}$ in order for the action of $\ee_7$ on a
single tower to close; these null elements of $\OO'$ are eigenstates of
$\eA_L$, the generator of the Cartan subalgebra of the $\sl(2,\RR)$ that
centralizes $\ee_7$ in $\ee_8$, as discussed further
in Section~\ref{summary}.

\subsection{Freudenthal product via commutators}
\label{fcomm}

It remains to show that this heuristic description is in fact correct.
We will need the multiplication table for null split octonions.  The elements
\begin{equation}
L_\pm = \frac{1\pm L}{2}
\end{equation}
satisfy
\begin{align}
L_\pm L_\pm &= L_\pm ,\\
L_\pm L_\mp &= 0.
\end{align}
and are thus projection operators.  We therefore introduce the sets
\begin{equation}
\NN_\pm
  = \{I_\pm,J_\pm,K_\pm\}
  = \{I,J,K\}L_\pm
\end{equation}
whose elements satisfy
\begin{align}
K_\pm K_\pm &= 0 ,\\
K_\pm K_\mp &= -L_\mp ,\\
K_\pm I_\pm &= J_\mp ,\\
K_\pm I_\mp &= 0 \label{KI0}
\end{align}
together with cyclic permutations of $\{I,J,K\}$.

We proceed in several steps, first considering orthogonal, tracefree elements
of the Albert algebra, then use Section~\ref{e8} to interpret the trace.
Let $\XXX\in\HHH$ be an element of the Albert algebra, and let
\begin{equation}
\XXX_0 = \XXX - \frac13\, (\tr\XXX)\,\III
\end{equation}
denote its tracefree part.  The Freudenthal product then decomposes as
\begin{align}
\XXX*\YYY
  &= \left(\XXX_0\circ\YYY_0-\frac12\,\tr(\XXX_0\circ\YYY_0)\,\III\right)
  \nonumber\\
  &\qquad
	-\frac16\Bigl( (\tr\XXX)\,\YYY_0 + (\tr\YYY)\,\XXX_0 \Bigr)
  \nonumber\\
  &\qquad
	+\frac19 (\tr\XXX)(\tr\YYY)\,\III
\end{align}
which implies that
\begin{align}
\XXX_0*\YYY_0 &= \XXX_0\circ\YYY_0-\frac12\,\tr(\XXX_0\circ\YYY_0)\,\III ,
 \label{s00}\\
x\,\III * y\,\III &= xy\,\III = x\,\III \circ y\,\III ,\label{s11}\\
x\,\III * \YYY_0 &= -\frac12\, x\,\YYY_0 \label{s10}.
\end{align}
A tracefree element $\XXX_0\in\HH$ maps directly to a \textit{labeled} Albert
algebra element $Q\XXX_0\in\ee_8$ for $Q\in\Im\OO'$, and we can ask what the
commutator of two such labeled elements is.  Since the trace of the Jordan
product is a positive-definite inner product on $\HHH$, there is a natural
orthonormal basis of (the tracefree elements of) $\HHH$, each of whose
elements squares to a diagonal matrix.  Working with such a basis, we can
assume that either $\XXX_0\perp\YYY_0$ or $\XXX_0=\YYY_0$.  In the first case,
direct computation establishes that
\begin{equation}
[Q_1\XXX_0,Q_2\YYY_0] = 2\,Q_1Q_2 (\XXX_0\circ\YYY_0)
\label{Mcomm}
\end{equation}
provided that $Q_1\perp Q_2$, so that $Q_2Q_1=-Q_1Q_2$, which we henceforth
assume.  Since we are assuming that the right-hand side of~\eqref{Mcomm} is
tracefree, it, too, is an element of $\ee_8$.  We have therefore shown that
\begin{equation}
\tr(\XXX_0\circ\YYY_0) = 0 \Longrightarrow
	[Q_1 \XXX_0,Q_2 \YYY_0] = 2\,Q_1Q_2 (\XXX_0*\YYY_0)
\label{XY00}
\end{equation}
where we have used~\eqref{s00} to replace the Jordan product by the
Freudenthal product.

In the second case, we can assume without loss of generality that
$\XXX_0=\YYY_0$ is of type I.  Explicitly, if
\begin{equation}
\EEE_I =
  \begin{pmatrix}z & e & 0 \\ \bar{e} & -z & 0 \\ 0 & 0 & 0 \\\end{pmatrix}
\end{equation}
is a basis element of type~I, then the normalization condition is
\begin{equation}
1 = \tr(\EEE_I\circ\EEE_I) = 2\,(z^2+|e|^2)
\end{equation}
and we have
\begin{align}
\EEE_I\circ\EEE_I
 &= \begin{pmatrix}\frac12 & 0 & 0 \\ 0 & \frac12 & 0 \\ 0 & 0 & 0 \\\end{pmatrix}
  = \frac13\,\III + \frac16
	\begin{pmatrix}1 & 0 & 0 \\ 0 & 1 & 0 \\ 0 & 0 & -2 \\\end{pmatrix}
\label{E1dot},\\
\EEE_I*\EEE_I
 &= \begin{pmatrix}0 & 0 & 0 \\0 & 0 & 0 \\ 0 & 0 & -\frac12 \\\end{pmatrix}
  = -\frac16\,\III + \frac16
	\begin{pmatrix}1 & 0 & 0 \\ 0 & 1 & 0 \\ 0 & 0 & -2 \\\end{pmatrix}
\label{E1star}
\end{align}
so that the matrix commutator becomes
\begin{equation}
[Q_1\EEE_I,Q_2\EEE_I]
  = 2\,Q_1Q_2 (\EEE_I\circ\EEE_I)
  = Q_1Q_2 (\III+2\,\EEE_I*\EEE_I)
\label{matcomm}
\end{equation}
with analogous results holding for basis elements of types II and III.

How are we to interpret this matrix commutator?  Section~\ref{e8} tells us the
answer: Decompose the result as the sum of a multiple of the identity matrix
and a tracefree matrix.  Thus, using~\eqref{E1dot} and~\eqref{matcomm},
\begin{equation}
[Q_1\EEE_I,Q_2\EEE_I]
  = \frac23\,Q_1Q_2 \,\III + \frac13\,\eS_{Q_1Q_2}
\end{equation}
which translates to a commutator in $\ee_8$ as
\begin{equation}
[Q_1\EEE_I,Q_2\EEE_I]
  = \frac23\,(Q_1\circ Q_2)\,\III + \frac13\,\eS_{Q_1Q_2} .
\end{equation}
\goodbreak
For example, with $Q_1=K_\pm=\frac12(K\pm KL)$ and $Q_2=I_\pm=\frac12(I\pm
IL)$, we have
\begin{align}
[K_\pm\,\EEE_I,I_\pm\,\EEE_I]
 &= \frac16\,\bigl((K\circ I+KL\circ IL)
    \nonumber\\
    &\qquad
	\pm(K\circ IL+KL\circ I)\bigr)\,\III
	+ \frac13\,\eS_{J_\mp}
\nonumber\\
 &= \frac16\,(L\circ JL\mp L\circ J)\,\III + \frac13\,\eS_{J_\mp}
\nonumber\\
 &= \frac{1}{6}\,\eG_{J_\mp} + \frac13\,\eS_{J_\mp}
\label{Xee}
\end{align}
where we have used the linear dependence of the 21 nested multiples of the
identity in the second equality, and then the relationship
\begin{align}
\eG_Q = -2(QL\circ L)\,\III
\label{Gid}
\end{align}
for $L\ne Q\in\Im\OO'$.  We have shown that
\begin{align}
[K_\pm\,\EEE_I,I_\pm\,\EEE_I]
 &= \pm\frac13\,(J_\mp\circ L)\,\III + \frac13\,\eS_{J_\mp} \nonumber\\
 &= -\frac13\,(J_\mp L\circ L)\,\III + \frac13\,\eS_{J_\mp}
\label{E1comm}
\end{align}
as well as similar results involving cyclic permutations of $I$, $J$, $K$.
Comparing~\eqref{E1comm} with~\eqref{E1star}, we can interpret the right-hand
side of~\eqref{E1comm} as ``$2\,J_\mp\,\EEE_I*\EEE_I$'' so long as we identify
$(QL\circ L)\,\III$ with $Q\,\III$, as previously proposed in~\cite{MagicE6}.
Thus, the right-hand side of~\eqref{XY00} holds for all trace-free
$\XXX_0,\YYY_0\in\HHH$.

What does this identification lead to for commutators involving the identity
matrix?  Using~\eqref{s11}, we clearly have
\begin{equation}
[Q_1\III,Q_2\III] = 2\,Q_1Q_2\,\III = 2\,Q_1Q_2\,\III*\III
\end{equation}
as matrices, assuming as usual that $Q_1\perp Q_2$, so we ask whether
\begin{equation}
[(Q_1L\circ L)\,\III,(Q_2L\circ L)\,\III] = 2\,(Q_1Q_2L\circ L)\,\III
\end{equation}
holds in~$\ee_8$ or, equivalently, whether
\begin{equation}
\left[\frac12\,\eG_{Q_1},\frac12\,\eG_{Q_2}\right] = -\eG_{Q_1Q_2}
\end{equation}
where we have used~\eqref{Gid}.  Although this relationship is not true for
general $Q_1$ and $Q_2$, it does hold if $Q_1,Q_2\in \NN_\pm$.  To establish
this result, it is sufficient to verify by direct computation that
\begin{equation}
[\eG_{K_\pm},\eG_{I_\pm}] = -4 \eG_{J_\mp}
\label{X11}
\end{equation}
and cyclic permutations.

What about commutators involving a single copy of the identity matrix, as
in~\eqref{s10}?  We must first interpret the matrix commutator
$[Q_1\III,Q_2\EEE_I]$
as the $\ee_8$ commutator
\begin{equation}
[(Q_1L\circ L)\,\III,Q_2\EEE_I]
  = \left[-\frac12\,\eG_{Q_1},Q_2\EEE_I\right]
\end{equation}
where we have again used~\eqref{Gid}.  We are again interested in the case
$Q_1,Q_2\in \NN_\pm$, and we compute directly that
\begin{equation}
\left[\eG_{K_\pm},I_\pm\,\EEE_I\right] = 2J_\mp\,\EEE_I
\label{X1e}
\end{equation}
and cyclic permutations.  We have shown that, in $\ee_8$,
\begin{equation}
[(Q_1L\circ L)\III,Q_2\EEE_I] = -Q_1Q_2\,\EEE_I = 2 Q_1Q_2\,\III*\EEE_I
\end{equation}
for $Q_1=K_\pm$ and $Q_2=I_\pm$, and hence for $Q_1,Q_2\in \NN_\pm$.
\footnote{This result differs from the matrix commutator
$[Q_1\III,Q_2\EEE_I]=2Q_1Q_2$, which can be attributed to the fact that
$(Q_1L\circ L)Q_2\ne Q_1Q_2$ as coefficients in $\ee_8$.}

We have therefore established the following lemma:
\begin{lemma}
Let $Q_1\ne Q_2\in \NN_\pm$ (with the same sign in both cases), and let
$\XXX,\YYY\in\HHH$ be elements of the Albert algebra.  Then
\begin{equation}
[Q_1\XXX,Q_2\YYY] = 2\,Q_1Q_2\,\XXX*\YYY
\end{equation}
in $\ee_8$, where the trace ``$Q\,\III$'' is to be interpreted as $(QL\circ
L)\,\III$, or equivalently as $-\frac12\eG_Q$.
\label{flemma}
\end{lemma}

\noindent
Lemma~\ref{flemma} can be generalized to the case $Q=Q_1=Q_2\in \NN_\pm$ by
noting that
\begin{align}
[Q\XXX_0,Q\YYY_0] &= 0 ,\\
[Q\XXX_0,\eG_Q] &= 0
\end{align}
and of course $Q^2=0$.

\subsection{Jordan trace via commutators}
\label{jcomm}

However, this is not the end of the story.  We can repeat the computations
leading up to~\eqref{XY00}, \eqref{Xee}, \eqref{X11}, and~\eqref{X1e}, but
with $I_\pm$ replaced by $I_\mp$.  Remembering~\eqref{KI0}, it should come as
no surprise that
\begin{align}
\tr(\XXX_0\circ\YYY_0) = 0 \Longrightarrow [K_\pm\,\XXX_0,I_\mp\,\YYY_0] &= 0 ,\\
\left[\eG_{K_\pm},I_\mp\,\EEE_I\right] &= 0 ,
\end{align}
but the analogs of~\eqref{Xee} and~\eqref{X11} require recomputation,
resulting in
\begin{align}
&[K_\pm\,\EEE_I,I_\mp\,\EEE_I]
\nonumber\\
 &= \frac16\,\bigl((K\circ I-KL\circ IL)
    \nonumber\\
    &\qquad
	\mp(K\circ IL-KL\circ I)\bigr)\,\III
	+ \frac13\,\eS_{K_\pm I_\mp}
\nonumber\\
 &= \frac{1}{2}\,\eA_{J_\pm} + 0 = \frac12\,\tr(\EEE_I\circ\EEE_I)\,\eA_{J_\pm} ,\\
&[\eG_{K_\pm},\eG_{I_\mp}]
 = 6 \eA_{J_\pm}
  = \frac12\,\tr(4\III)\,\eA_{J_\pm} .
\end{align}
%
Since $\tr(\EEE_I)=0$, we have established:
\begin{lemma}
Let $\XXX,\YYY\in\HHH$ be elements of the Albert algebra.  Then in $\ee_8$ we
have
\begin{equation}
[K_\pm\,\XXX,I_\mp\,\YYY] = \frac12\,\tr(\XXX\circ\YYY)\,\eA_{J_\pm}
\end{equation}
and its cyclic permutations in $I$, $J$, $K$, where the trace ``$Q\,\III$'' is
to be interpreted as $(QL\circ L)\,\III$, or equivalently as $-\frac12\eG_Q$.
\label{trlemma}
\end{lemma}

A straightforward consequence of Lemmas~\ref{flemma} and~\ref{trlemma} is that
the determinant of an element of the Albert algebra can be represented using
commutators in $\ee_8$.

\begin{corollary}
Let $\XXX\in\HHH$ be an element of the Albert algebra.  Then in $\ee_8$ we
have
\begin{equation}
\bigl[K_\pm\XXX,[K_\pm\XXX,I_\pm\XXX]\bigr] = -3\,(\det\XXX)\,\eA_{I_\mp}
\end{equation}
where the trace ``$Q\,\III$'' is to be interpreted as $(QL\circ L)\,\III$, or
equivalently as $-\frac12\eG_Q$.
\label{CDet}
\end{corollary}

\begin{proof}
Using Lemmas~\ref{flemma} and~\ref{trlemma},
\begin{align}
\bigl[K_\pm\XXX,[K_\pm\XXX,I_\pm\XXX]\bigr]
  &= \bigl[K_\pm\XXX,2J_\mp(\XXX*\XXX)\bigr] \nonumber\\
  &= -2\bigl[J_\mp(\XXX*\XXX),K_\pm\XXX\bigr] \nonumber\\
  &= -\tr\bigl((\XXX*\XXX)\circ\XXX\bigr)\,\eA_{I_\mp} \nonumber\\
  &= -3\,(\det\XXX)\,\eA_{I_\mp}
\label{detcomm}
\end{align}
where the last step uses the definition of the determinant in the Albert
algebra.
\label{cor1}
\end{proof}

\subsection{Anchors}
\label{anchors}

Lemmas~\ref{flemma} and~\ref{trlemma} are the heart of our construction of
(two copies of) the minimal representation of $\ee_7$.  These lemmas show how
to interpret the elements $K_\pm\XXX\in\ee_7$ as raising and lowering
operators, acting on copies of the Albert algebra, namely $I_\pm\XXX$ and
$J_\pm\XXX$.  Comparison with~\eqref{Araise} and~\eqref{Braise} justifies the
identification of $K_\pm\XXX$ with the null translations of $\ee_7$, as
claimed in Section~\ref{overview}.

Comparing Lemma~\ref{trlemma} with Table~\ref{tower} then also shows that the
anchors of our two Freudenthal towers must be (multiples of) $\eA_{I_\pm}$ and
$\eA_{J_\pm}$.  It remains to verify that the null translations $K_\pm\XXX$ act
correctly on the anchors, that is, in agreement with the first maps
in~\eqref{Araise} and~\eqref{Braise}.

We begin by investigating $[K_\pm\XXX,\eA_{J_\pm}]$, where the signs are the
same.  Using Lemma~\ref{trlemma} and the Jacobi identity (and of course
interpreting $Q\,\III$ as $-\frac12\eG_Q$), we have
\begin{align}
[K_\pm\XXX,\eA_{J_\pm}]
 &= \left[K_\pm\XXX,\frac23\,[K_\pm\III,I_\mp\III]\right] \nonumber\\
 & = \frac23\, \bigl[[K_\pm\XXX,K_\pm\III],I_\mp\III\bigr]
	\nonumber\\
	&\qquad
	+ \frac23\, \bigl[K_\pm\III,[K_\pm\XXX,I_\mp\III]\bigr] \nonumber\\
 &= 0 + \left[K_\pm\III,\frac13\,(\tr\XXX)\,\eA_{J_\pm}\right]
\label{KpAJp1}
\end{align}
since the first inner commutator vanishes using (the generalization of)
Lemma~\ref{flemma} and the fact that $K_\pm^2 = 0$, and the second can be
evaluated using Lemma~\ref{trlemma}.  Finally, it can be verified that
\begin{equation}
[K_\pm\III,\eA_{J_\pm}] = 0
\label{KpAJp2}
\end{equation}
by direct computation.

Turning to $[K_\pm\XXX,\eA_{J_\mp}]$, where the signs are now different, and
proceeding as in~\eqref{KpAJp1}, we obtain
\begin{align}
[K_\pm\XXX,\eA_{J_\mp}]
 &= \left[K_\pm\XXX,\frac23\,[K_\mp\III,I_\pm\III]\right] \nonumber\\
 &= \frac23\,\bigl[[K_\pm\XXX,K_\mp\III],I_\pm\III\bigr]
	\nonumber\\
	&\qquad
	+ \frac23\,\bigl[K_\mp\III,[K_\pm\XXX,I_\pm\III]\bigr] .
\end{align}
The last term can be evaluated using Lemma~\ref{flemma}, since
\begin{align}
\bigl[K_\mp\III,[K_\pm\XXX,I_\pm\III]\bigr]
 &= [K_\mp\III,2J_\mp\>\III*\XXX] \nonumber\\
 &= -4I_\pm\>\III*(\III*\XXX) .
\end{align}
However, the first term requires new results, namely that
\begin{equation}
[K_\pm\XXX,K_\mp\III] = \pm L\XXX
\label{new1}
\end{equation}
and
\begin{equation}
[L\XXX,I_\pm\III] = \mp I_\pm (2\XXX-\frac13(\tr\XXX)\,\III)
\label{new2}
\end{equation}
which follow either from straightforward generalizations of Lemmas~\ref{flemma}
and~\ref{trlemma} or by direct computation.
\footnote{The definition of $\eG_L$, and hence of ``$L\,\III$'', is
conventional, and should more accurately be written in~\eqref{new1}
and~\eqref{new2} as $-2(K\circ KL)\,\III$.  Thus, the interpretation of the
trace term in ``$L\,\XXX$'' depends on the labels (here $K$ and $KL$) being
considered.  With this caveat, cyclic permutations of~\eqref{new1}
and~\eqref{new2} also hold.}
Since
\begin{equation}
\III*(\III*\XXX) = \frac14\bigl(\XXX+(\tr\XXX)\,\III\bigr)
\end{equation}
we conclude that
\begin{align}
[K_\pm\XXX,\eA_{J_\mp}]
  &= \frac23\,\bigl(-(\XXX+(\tr\XXX)\,\III)
	\nonumber\\
	&\qquad
	-(2\XXX-(\tr\XXX)\,\III)\bigr) I_\pm
\nonumber\\
  &= -2I_\pm\XXX 
\label{KpAJm}
\end{align}
in agreement with~\eqref{Araise} and~\eqref{Braise}.

\begin{table}
\centering
\begin{tabular}{|c|c|c|}
\hline
& \boldmath$\eA_L$ & \boldmath$\eG_L$ \\
\hline
\boldmath$I+IL$ & $-2$ & $-2$ \\
\boldmath$J+JL$ & $+2$ & $-2$ \\
\boldmath$K+KL$ & $~~\,0$ & $+4$ \\
\boldmath$I-IL$ & $+2$ & $+2$ \\
\boldmath$J-JL$ & $-2$ & $+2$ \\
\boldmath$K-KL$ & $~~\,0$ & $-4$ \\
\boldmath$L$ & $~~\,0$ & $~~\,0$ \\
\hline
\end{tabular}
\caption{The eigenvectors and corresponding eigenvalues of the generators
$\{\eA_L,\eG_L\}$ of the Cartan subalgebra of $\gg_{2(2)}$ when acting on
$\Im\OO'$.}
\label{g2ev}
\end{table}

For completeness, we note explicitly that interchanging $I$ and $J$ results in
\begin{align}
[K_\pm\XXX,A_{I_\pm}] &= +2J_\pm\XXX ,\label{KpAIp}\\
[K_\pm\XXX,A_{I_\mp}] &= 0 ,\label{KpAIm}
\end{align}
by a similar argument.

Examining~\eqref{KpAJp1}--\eqref{KpAJp2}, \eqref{KpAJm}, and
\eqref{KpAIp}--\eqref{KpAIm}, together with Lemmas~\ref{flemma}
and~\ref{trlemma}, and comparing with~\eqref{FreudX}--\eqref{qFreud}, it is
straightforward to conclude that the anchors are given by $\pm\frac12 A_Q$,
with $Q\in\{I_\pm,J_\pm\}$.  However, getting the signs to match exactly with
those in~\eqref{FreudX}--\eqref{qFreud} requires some care.  One possibility
is to choose the anchors to be $-\frac12 A_Q$, then use
$\{\pm K_\pm\AAA,\mp K_\mp\BBB\}$ as the raising and lowering operators for
the towers with anchors labeled by $\{I_\mp,J_\mp\}$, respectively, as shown
schematically in Table~\ref{e8tower}.

\subsection{The dilation}
\label{dilation}

As described in Section~\ref{e7Freud}, $\ee_7$ consists of $\ee_6$ together
with two sets of null translations, as well as a single dilation.  The action
of $\ee_6$ on the rest of adjoint $\ee_8$ was discussed in~\cite{MagicE6}, and
we showed how null translations act in Sections~\ref{fcomm} and~\ref{jcomm}.
It remains only to consider the action of the dilation on Freudenthal towers.

However, the dilation corresponds to (a multiple of) the element $\eG_L$,
which is in our Cartan subalgebra, and which has the same eigenvectors as
$\eA_L\in\sl(2,\RR)$.  More precisely, $\{\eA_L,\eG_L\}$ generates the
Cartan subalgebra of $\gg_{2(2)}$, whose minimal representation is just
$\Im\OO'$, and whose common eigenvectors are precisely our labels $\NN_\pm$.
Our conventions have been chosen so that the eigenvalues are as given in
Table~\ref{g2ev}.  The action of the dilation on any element of $\ee_8$ of the
form $Q\XXX$, with $Q\in\Im\OO'$ and $\XXX\in\HHH$, is given by the
$\gg_{2(2)}$ action on~$Q$; if $Q\in \NN_\pm$, $Q\XXX$ is therefore an
eigenvector of $\eG_L$.
\footnote{This result is obvious if $\tr\XXX=0$.  If not, we make the usual
identification of $Q\III$ with $-\frac12\eG_Q$, and we verify by direct
computation that
\[
[\eA_L,\eG_Q] = \lambda \eG_Q
\]
for $Q\in \NN_\pm$ and $\lambda$ chosen appropriately from the middle column
of Table~\ref{g2ev}.}

However, Freudenthal towers also contain elements that are not of the form
$Q\XXX$, namely the anchors ($A_{I_\pm}$ and $A_{J_\pm}$).  We compute
explicitly that
\begin{align}
[\eG_L,\eA_{I_\pm}] &= \mp6 \eA_{I_\pm} ,\label{GAI}\\
[\eG_L,\eA_{J_\pm}] &= \pm6 \eA_{J_\pm} .\label{GAJ}
\end{align}
Comparing Table~\ref{g2ev} and~\eqref{GAI}--\eqref{GAJ}
with~\eqref{FreudX}--\eqref{qFreud}, we see that the correctly scaled dilation
is given by $-\frac16\eG_L$.

\subsection{The action of $\ee_6$}
\label{e6}

The careful reader will have noticed that we have not yet discussed the action
of $\ee_6$ ($\phi$) in~\eqref{ThetaDef}.  However, this action is
straightforward.

First of all, $\ee_6$ commutes with the anchors, which are elements of
$\sl(3,\RR)$, the centralizer of $\ee_{6(-26)}$ in $\ee_{8(-24)}$.  This
property of $\ee_8$ correctly reflects the fact that $\phi$ only acts on
$\XXX$ and $\YYY$ in~\eqref{FreudX}--\eqref{qFreud}, and not on $p$ or $q$.

The action of $\ee_6$ on labeled Albert algebras via commutators in $\ee_8$
was discussed in detail in~\cite{MagicE6}.  In a nutshell, we can consider
elements $\phi\in\ee_6$ to be $3\times3$ matrices over $\OO$, with rotations
corresponding to tracefree anti-Hermitian matrices and boosts to tracefree
Hermitian matrices; the derivations in $\gg_2\subset\ee_6$ can be handled as
nested rotations.  In both cases, $\phi$ acts on $\XXX\in\HHH$ via
$\XXX\longmapsto\phi\XXX+\XXX\phi^\dagger$.  In $\ee_8$, $\XXX$ acquires a
label such as $K_\pm$.  If $\phi$ is a rotation, then we can identify $\phi$
with the element $\phi_0=\phi\in\ee_8$, since the action is already given by
the commutator, and the label pulls through the calculation.  Explicitly, we
have
\begin{equation}
[\phi_0,K_\pm\XXX] = K_\pm[\phi,\XXX] = K_\pm \phi(\XXX) .
\end{equation}
However, if $\phi$ is a boost, we must identify instead
$\phi_0=L\phi\in\ee_8$, since
\begin{equation}
[L\phi,K_\pm\XXX]
  = \mp2 K_\pm \phi\circ\XXX
  = \mp K_\pm \phi(\XXX) .
\end{equation}
Since $\phi$ and $\phi'$ differ precisely in the sign of the action of boosts,
we have shown that
\begin{equation}
[\phi_0,K_-\XXX] = K_-\phi(\XXX) ,
\quad
[\phi_0,K_+\YYY] = K_+\phi'(\YYY) ,
\end{equation}
in agreement with~\eqref{FreudX} and~\eqref{FreudY}.

\subsection{\boldmath Summary of the action on Freudenthal towers}
\label{summary}


\begin{table}
\centering
\fbox{
\begin{tabular}{cccc}
& $p\eA_{I_+}$ & $p\eA_{J_-}$ \\
& $|$ & $|$ \\
$\uparrow$\quad\null & $J_+\YYY$ & $I_+\YYY$ \\
$K_-\AAA$\quad\null & $|$ & $|$ & \quad$K_+\BBB$ \\
& $I_-\XXX$ & $J_-\XXX$ & \quad$\downarrow$ \\
& $|$ & $|$ \\
& $q\eA_{J_+}$ & $q\eA_{I_-}$
\end{tabular}
}
\caption{A schematic representation of the two Freudenthal towers contained in
adjoint~$\ee_8$, showing also the action of the null translations in $\ee_7$.
Each tower can be mapped to \hbox{$\Pc=(\XXX,\YYY,p,q)$}, and the null
translations to $(0,0,\AAA,0)$ and $(0,0,0,\BBB)$, as in
Section~\ref{e7Freud}.}
\label{e8tower}
\end{table}

Recall that we have chosen our $\ee_7$ subalgebra such that the Albert
algebras labeled by $K_\pm$ are in $\ee_7\subset\ee_8$, whereas those labeled
by $I_\pm$ and $J_\pm$ are not, nor are $\eA_{I_\pm}$ or $\eA_{J_\pm}$.  We
have shown that $\ee_7$ acts on these $4\times27+4=112$ elements as two
Freudenthal towers, as shown schematically in Table~\ref{e8tower}.
Explicitly, we have identified Freudenthal's description~\eqref{ThetaDef} of an
element $\Theta$ of $\ee_7$ with the $\ee_8$ element
\begin{equation}
 \Theta_0 = \phi_0 + \frac16\,\rho\,\eG_L + K_-\AAA - K_+\BBB
\label{FreudTheta}
\end{equation}
while also identifying Freudenthal's description~\eqref{Pcdef} of an element
$\PPP$ of the minimal representation of $\ee_7$ with either of the $\ee_8$
elements
\begin{align}
 \PPP_+ &= I_-\XXX + J_+\YYY - \frac12\,p\eA_{I_+} - \frac12\,q\eA_{J_+}
	,\label{FreudPp}\\
 \PPP_- &= J_-\XXX - I_+\YYY - \frac12\,p\eA_{J_-} + \frac12\,q\eA_{I_-}
	\label{FreudPm} .
\end{align}



%

The centralizer of $\ee_{7(-25)}$ in $\ee_{8(-24)}$ is $\sl(2,\RR)$, generated
by $\{\eA_K,\eA_{KL},\eA_L\}$.  Choosing $\eA_L$ to generate its Cartan
subalgebra, the raising and lower operators are the null combinations
$\eA_{K\mp  KL}$, respectively.  Our construction of Freudenthal towers
is adapted to this description of $\sl(2,\RR)$, since
\begin{equation}
[\eA_L,\PPP_\pm] = \pm2 \PPP_\pm .
\end{equation}
Thus, for fixed $\PPP$, the pair $\{\PPP_\pm\}$ generates a (real) doublet of
$\sl(2,\RR)$.

\section{\boldmath The structure of $\ee_7$}
\label{e7adj}


We return to the structure of adjoint $\ee_7\subset\ee_8$.
Freudenthal~\cite{FreudenthalE7} tells us that the commutator of two elements
$\Theta_1,\Theta_2\in\ee_7$ of the form~\eqref{ThetaDef} is another such
element $\Theta$, with 
\begin{align}
  \phi
	&= [\phi_1,\phi_2] - 2\langle\AAA_1,\BBB_2\rangle
		+ 2\langle\AAA_2,\BBB_1\rangle
       ,\label{FreudPhi}\\
  \rho &= -\tr(\AAA_1\circ\BBB_2) + \tr(\AAA_2\circ\BBB_1) ,\label{FreudRho}\\
  \AAA &= \left(\phi_1-\frac23\rho_1\right)\AAA_2
         - \left(\phi_2-\frac23\rho_2\right)\AAA_1 ,\label{FreudA}\\
  \BBB &= \left(\phi'_1+\frac23\rho_1\right)\BBB_2
         - \left(\phi'_2+\frac23\rho_2\right)\BBB_1 ,\label{FreudB}
\end{align}
where $\phi'$ denotes the dual of $\phi$, and the $\ee_6$ element
$\langle\AAA,\BBB\rangle$ acts on the Albert algebra via
\begin{align}
\langle\AAA,\BBB\rangle:
  \XXX &\longmapsto \BBB\circ(\AAA\circ\XXX) - \AAA\circ(\BBB\circ\XXX)
  \nonumber\\
  &\qquad - (\AAA\circ\BBB)\circ\XXX + \frac13\tr(\AAA\circ\BBB)\XXX .
\label{e6AB}
\end{align}
When using~\eqref{FreudX}--\eqref{qFreud} to verify~\eqref{FreudPhi}, the
alternate form
\begin{align}
\langle\AAA,\BBB\rangle:
  \XXX &\longmapsto 2\bigl(\BBB*(\AAA*\XXX)\bigr)
	- \frac12 \tr(\BBB\circ\XXX)\,\AAA
  \nonumber\\
  &\qquad - \frac16\tr(\AAA\circ\BBB)\XXX
\end{align}
is useful~\cite{FreudenthalE7}.  It is also worth observing that if
$\BBB=\AAA*\AAA$ then the first two terms in~\eqref{e6AB} cancel, as do the
last two terms, so that
\begin{equation}
\langle\AAA,\AAA*\AAA\rangle=0 .
\label{AstarA}
\end{equation}

We emphasize that the commutator of two elements of $\ee_7\subset\ee_8$ of the
form~\eqref{FreudTheta} \textit{must}
reproduce~\eqref{FreudPhi}--\eqref{FreudB}, since we have verified that their
action on $\PPP_\pm\in\ee_8$ reproduces~\eqref{FreudX}--\eqref{qFreud}.  We
can nonetheless verify by direct computation in $\ee_8$
that~\eqref{FreudPhi}--\eqref{FreudB} hold.  Most of this computation is
straightforward, since each term in~\eqref{FreudTheta} is separately an
eigenvector of $\eG_L$, and the action of $\ee_6$ ($\phi_0$) on labeled
Albert algebras was summarized in Section~\ref{e6}.

The final step is to compute $[K_-\AAA,K_+\BBB]$.  Naively, since
$K_\pm K_\mp=-L_\mp$, we expect to have something like
$[K_-\AAA,K_+\BBB]
  = -\frac12[\AAA,\BBB] - L(\AAA\circ\BBB)
$.
However, due to the lack of associativity, the correct expansion turns out to
be
\begin{equation}
[K_-\AAA,K_+\BBB]
  = -\frac12[L\AAA,L\BBB] - L(\AAA\circ\BBB) .
\end{equation}
Lengthy but straightforward computation now shows that if the label $Q$
anticommutes with $L$, then
\begin{align}
  \bigl[[L\AAA&,L\BBB],Q\XXX\bigr]
   \nonumber\\
  &= \bigl[[L\AAA,Q\XXX],L\BBB\bigr] + \bigl[L\AAA,[L\BBB,Q\XXX\bigr]]
   \nonumber\\
  &= (LQ\,\AAA\XXX-QL\,\XXX\AAA)L\BBB - L\BBB(LQ\,\AAA\XXX-QL\,\XXX\AAA)
   \nonumber\\
  &\qquad + L\AAA(LQ\,\BBB\XXX-QL\,\XXX\BBB)
   \nonumber\\
  &\qquad - (LQ\,\BBB\XXX-QL\,\XXX\BBB)L\AAA
   \nonumber\\
  &= 2\bigl( -(\AAA\circ\XXX)\,QL\,L\BBB - L\,LQ\,\BBB\,(\AAA\circ\XXX)
   \nonumber\\
  &\qquad + L\,LQ\,\AAA\,(\BBB\circ\XXX) + (\BBB\circ\XXX)\,QL\,L\AAA \bigr)
   \nonumber\\
  &= 2\bigl( -(\AAA\circ\XXX)\,Q\,\BBB - \,Q\,\BBB\,(\AAA\circ\XXX)
   \nonumber\\
  &\qquad + \,Q\,\AAA\,(\BBB\circ\XXX) + (\BBB\circ\XXX)\,Q\,\AAA \bigr)
   \nonumber\\
  &= 4Q\,\bigl( -(\AAA\circ\XXX)\circ\BBB + \AAA\circ(\BBB\circ\XXX) \bigr)
   \nonumber\\
  &= 4Q\,\bigl( \AAA\circ(\BBB\circ\XXX) - \BBB\circ(\AAA\circ\XXX) \bigr)
\end{align}
and a similar computation shows that
\begin{align}
\Bigl[2L(\AAA\circ\BBB) &+ \frac13\tr(\AAA\circ\BBB)\eG_L,Q\XXX\Bigr]
   \nonumber\\
  &= \pm 4Q\Bigl(
	\bigl(\AAA\circ\BBB - \frac13\tr(\AAA\circ\BBB)\III\bigr)\circ\XXX
	\Bigr)
\end{align}
with the additional assumption that $LQ=\pm Q$.  Combining these results, we
obtain, for instance with $LQ=Q$,
\begin{align}
  \bigl[[K_-\AAA&,K_+\BBB],I_-\XXX\bigr]
   \nonumber\\
  &= 2I_-\langle\AAA,\BBB\rangle(\XXX)
     +\frac16\tr(\AAA\circ\BBB)[\eG_L,I_-\XXX]
\label{Kmp}
\end{align}
with the same relation holding if $I_-$ is replaced by~$J_-$.
These results verify both~\eqref{FreudPhi} and~\eqref{FreudRho}, taking into
account the minus sign in~\eqref{FreudTheta}.  (Using labels such as $Q=I_+$
satisfying $LQ=-Q$ yields a similar result, but with the dual action of
$\langle\AAA,\BBB\rangle$.  Our conventions are chosen here so as to match
those of Freudenthal~\cite{FreudenthalE7}.)

\section{Discussion}
\label{discussion}

\subsection{Graded Structure}
\label{graded}

Our decomposition of $\ee_8$ into Freudenthal towers over
$\ee_7\oplus\sl(2,\RR)$ is closely related to the graded structures discussed
by Bars and G\"unaydin~\cite{BarsGunaydin79,Gunaydin89}, building on previous
work by Kantor~\cite{Kantor72,Kantor73}.  As noted by G\"unaydin and
Hyun~\cite{Gunaydin88}, $\ee_{7(-25)}$ admits a graded structure of the form
\begin{equation}
\gg = ... \gg_{-1} \oplus \gg_0 \oplus \gg_1 ...
\end{equation}
with $\gg_0$ semisimple, $\gg_m$ nilpotent for $m\ne0$, and
\begin{equation}
[\gg_m,\gg_n] \subset \gg_{m+n}
\end{equation}
namely
\begin{equation}
\ee_{7(-25)}
  = \bar{\mb{27}} \oplus \bigl(\ee_{6(-26)}\oplus\so(1,1)\bigr) \oplus \mb{27} .
\end{equation}
As can be verified by direct but somewhat messy computation, the Lie algebra
$\ee_{8(-24)}$ admits a similar graded structure over $\ee_{7(-25)}$, namely
\begin{equation}
\ee_{8(-24)}
  = \bar{\mb{56}} \oplus \bigl(\ee_{7(-25)}\oplus\sl(2,\RR)\bigr)
	\oplus \mb{56},
\end{equation}
as well as a ``dual'' grading over $\ee_{6(-26)}$, namely
\begin{align}
\ee_{8(-24)}
  &= (2\times\mb{1}) \oplus \mb{27} \oplus (2\times\bar{\mb{27}})
\nonumber\\
  & \qquad \oplus \bigl(\ee_{6(-26)}\oplus\so(1,1)\oplus\sl(2,\RR)\bigr)
\nonumber\\
  & \qquad \oplus (2\times\mb{27}) \oplus \bar{\mb{27}} \oplus (2\times\mb{1})
\end{align}
in which both $\ee_{7(-25)}$ and the two Freudenthal towers
($\bar{\mb{56}}\oplus\mb{56}$) are chopped up into their component pieces and
rearranged.

\subsection{Other Signatures}
\label{signature}

The construction of $\ee_8$ given in Section~\ref{e8} works for all three real
forms of $\ee_8$, namely
\begin{align}
\ee_8 = \ee_{8(-248)} &\isom \su(3,\OO\otimes\OO) ,\\
\ee_{8(-24)} &\isom \su(3,\OO'\otimes\OO) ,\\
\ee_{8(8)} &\isom \su(3,\OO'\otimes\OO') ,
\end{align}
depending on the choice of division algebras used.
Three real forms of $\ee_7$ are contained in
$\ee_8\isom\su(3,\OO'\otimes\OO)$, namely
$\ee_{7(-25)}\isom\su(3,\HH'\otimes\OO)$,
$\ee_{7(-5)}\isom\su(3,\OO'\otimes\HH)$, and compact
$\ee_7\isom\su(3,\HH\otimes\OO)$; the split real form
$\ee_{7(7)}\isom\su(3,\OO'\otimes\HH')$ is contained in the split real form
$\ee_{8(8)}\isom\su(3,\OO'\otimes\OO')$.

However, the representation theory is slightly different in the cases
involving $\HH$ rather than $\HH'$, where we expect the decomposition of
$\ee_8$ over $\ee_7\oplus\su(2)$ to be
\begin{equation}
\mb{248} = \mb{133} + 1\times\mb{112} + \mb{3}
\end{equation}
with one \textit{complex} Freudenthal tower, rather than two \textit{real}
Freudenthal towers.  This behavior is easy to track in our notation, since the
split into real Freudenthal towers relied on the existence of projection
operators $1\pm L$, which only exist in $\OO'$.

\subsection{The Standard Model of Particle Physics}
\label{SM}

In~\cite{Octions}, we proposed interpreting $\sl(3,\RR)$ as the color symmetry
in the Standard Model.  As noted in Section~\ref{e7}, an intriguing
consequence of~\eqref{e6AB} is~\eqref{AstarA}, which, when inserted
in~\eqref{Kmp}, yields another expression for the determinant, namely
\begin{equation}
[K_-\AAA,K_+(\AAA*\AAA)] = \frac12(\det\AAA)\,\eG_L  .
\label{Kdet}
\end{equation}
Furthermore, Lemma~\ref{flemma} says that
\begin{equation}
2K_+(\AAA*\AAA) = [I_-\AAA,J_-\AAA] ,
\end{equation}
so that~\eqref{Kdet} becomes
\begin{equation}
\bigl[K_-\AAA,[I_-\AAA,J_-\AAA]\bigr] = (\det\AAA)\,\eG_L .
\end{equation}
which can be polarized, yielding a colorless triple product on the Albert
algebra.  As discussed in~\cite{Other}, every Albert algebra element can be
uniquely decomposed into primitive idempotents, playing the role of
eigenvectors---and the triple product of these idempotents reproduces the
determinant of the original matrix.  Do these properties provide a mechanism
for constructing colorless baryons?  If so, does~\eqref{Kdet} itself provide a
similar mechanism for constructing colorless mesons?

\subsection{Conclusion}
\label{conclude}

It is clear from representation theory and the decomposition~(\ref{decomp})
that the induced action using commutators in $\ee_8$ of $\ee_7$ on each of the
$\mb{56}$s in the decomposition must reproduce the action of $\ee_7$ on
Freudenthal towers.  That action, summarized in Section~\ref{e7Freud},
involves both the Jordan and Freudenthal products on the Albert algebras
contained in the $\mb{56}$s, as we verified explicitly in Sections~\ref{jcomm}
and~\ref{fcomm}.  Although this equivalence was guaranteed on theoretical
grounds, it is nonetheless remarkable to see explicitly how commutators in
$\ee_8$ reproduce (much of) the algebraic structure of the Albert algebra.

\newpage
\bibliographystyle{unsrt}
\bibliography{octo,octo2,octo3}

\end{document}